\documentclass[xypic,amscd,syntonly,amssymb]{amsart}

\usepackage{graphicx}
\usepackage[all]{xy}
\usepackage{amssymb}
\usepackage{amsmath}
\usepackage{epsfig}
\theoremstyle{plain}
\newtheorem{Thm}{Theorem}
\newtheorem{Prop}[Thm]{Proposition}

\newtheorem{Lem}[Thm]{Lemma}

 \theoremstyle{definition}
\theoremstyle{remark}

\numberwithin{equation}{section}

\begin{document}
 \title{ Action Integrals and Infinitesimal Characters}

 \author{ ANDR\'{E}S   VI\~{N}A}
\address{Departamento de F\'{i}sica. Universidad de Oviedo.   Avda Calvo
 Sotelo.     33007 Oviedo. Spain. }
 \email{vina@uniovi.es}
\thanks{This work has been partially supported by Ministerio de Educaci\'on y
 Ciencia, grant MAT2007-65097-C02-02}

 \maketitle
\begin{abstract}
Let $G$ be a reductive Lie group and ${\mathcal O}$ the coadjoint
orbit of a hyperbolic element of ${\frak g}^*$. By $\pi$ is
denoted the unitary irreducible representation of $G$ associated
with ${\mathcal O}$  by the orbit method. We give geometric
interpretations in terms of concepts related to ${\mathcal O}$ of
the constant $\pi(g)$, for $g\in Z(G)$.
 We also offer a description of the invariant $\pi(g)$ in terms of
 action integrals and Berry phases.
 In the spirit of the orbit method we interpret geometrically the
infinitesimal character of the differential representation of
$\pi$.

\end{abstract}
   \smallskip

%\subjclass[2000]{ 53D50, 22E45}


%%%%%%%%%%%%%%%%%%%%%%%%%%%%%%%%%%%%%%%%%%%%%%%%%%%%%%%%%%%%%%%%%%%%%%%%%%%%%%%%%%%%%%%%%%%
%%%%%%%%%%%%%%%%%%%%%%%%%%%%%%%%%%%%%%%%%%%%%%%%%%%%%%%%%%%%%%%%%%%%%%%%%%%%%%%%%%%%%%%%%%%%%%%%%%%%
\section{Introduction}
Roughly speaking, the orbit method \cite{Kir}, \cite{Vo2} suggests
that the unitary dual of a Lie group $G$ (i.e. the set of
equivalence classes of unitary irreducible representations of $G$)
is in bijective correspondence with the space of coadjoint orbits
of $G$.
 %pairs consisting
 %of a coadjoint orbit of $G$ and an ``integral datum" (see Section
 %\ref{S:definitions}).
 Moreover the orbit method relates geometric
properties of the coadjoint orbit with properties of the
corresponding irreducible representation. This bijective
correspondence exists if $G$ is a connected simply connected
nilpotent group; in other cases where the correspondence is not a
perfect bijection this method gives valuable suggestions about the
geometric meaning  of some facts of representation theory.

In this paper $G$ will be a reductive group and ${\mathcal O}$
will be a coadjoint orbit of  a hyperbolic element $\eta\in{\frak
g}^*$, where ${\frak g}^*$ is the dual of the Lie algebra ${\frak
g}$ of $G$. In the spirit of the orbit method we will
 give geometric interpretations of some invariants of the representation associated with
 ${\mathcal O}$. This will allow us, in turn, to offer physical
 interpretations of those invariants in terms of action integrals
 and Berry phases along curves generated in  physical systems by the action of symmetry groups.
 This is valid for   groups relevant in Physics, such as:
 $SO(p,q)$, $Sp(2n)$, $SL(n,{\mathbb R})$, etc.

\smallskip

   For the construction of a representation of $G$ from the
 orbit
 %it is necessary to impose an integrality condition to the
 %orbit; in the hyperbolic case
 we will assume that ${\mathcal O}$
 admits an integral datum (see Section
 \ref{S:definitions}). By
 means of
 an integral datum    one defines
  a unitary
irreducible representation $\pi$ of $G$ by induction from a
parabolic subgroup of $G$.
%We will   consider the following points:
 According to Schur's lemma, if $g_1$ belongs to the center of
$G$, $\pi(g_1)$ is a scalar operator defined by a constant
$\kappa$,
\begin{equation}\label{kappadefinition}
\pi(g_1)=
 \kappa\,{\rm Id}.
  \end{equation}
We will give
%a geometric meaning to the constant $\kappa$.
%In the
%spirit of the orbit method we will
   %express in geometric terms:
  an interpretation in geometric terms of:

   1) The constant $\kappa$.

    2) The infinitesimal character
    \cite{Kn} of $\pi'$, the differential
representation of $\pi$, considered as a representation of the
universal enveloping algebra $U({\frak g})$ of ${\frak g}$.

 3) Some values of the character $\chi_{\tau}$ of $\tau$, where $\tau$ is any irreducible
 representation of any maximal compact
 subgroup $K$  of $G$, which occurs in $\pi_{|K}$.
  %The values of the character $\chi_{\pi}$ of $\pi$ at elements of the
 %form $ala^{-1}$, with $l$ in the stabilizer of $\eta$ for the
 %coadjoint action and $a\in G$.

The orbit ${\mathcal O}$ is the homogeneous space $G/L$ , where
$L$ is the stabilizer of $\eta$.  An integral datum is a unitary
irreducible representation of $L$ on a Hilbert space $H$,
satisfying an additional condition (see (\ref{datum})). Given an
integral datum $\Lambda$, by $\Phi$ we denote the representation
of $L$, tensor product of $\Lambda$ and the character   of $L$
 on half-densities (\ref{Delta}).

 As $\eta$ is a hyperbolic element
its orbit ${\mathcal O}$ possesses a real polarization defined by
a subalgebra ${\frak u}$ of ${\frak g}$ (see (\ref{fraku})).
 ${\mathcal B}_1$ will be the space of smooth $\Phi$-equivariant maps $s:G\to
 H$, with compact support modulo $L$, and such that $L_C(s)=0$,
 for $C\in{\frak u}$, where $L_A$ is the left invariant vector
 field on $G$ determined by $A \in{\frak g}$.
 The representation $\pi$ is  the left regular
 representation of $G$ defined on
  the completion of the pre-Hilbert space ${\mathcal
 B}_1$.  Thus  the operator associated to $A \in{\frak g}$ in the  differential representation  $\pi'$ of
 $\pi$ is $-R_A$, the opposite of the right invariant vector field
 on $G$ defined by $A$.

 As a first step we will define the representation $\pi'$ in the
 context of  fibre bundles.
We will consider the $GL(H)$-principal bundle ${\mathcal
F}:=G\times_{\Phi} GL(H)$ over ${\mathcal O}$,
 defined by means of the representation $\Phi$ of $L$.
On ${\mathcal F}$ there is a natural left $G$-action
 and an obvious ${\mathbb C}^*$-action induced by the multiplication by  nonzero
scalars of elements of $GL(H)$.
 In particular each $A\in{\frak g}$
defines a vector field $Y_A$ on ${\mathcal F}$ by the $G$-action.
Furthermore, on ${\mathcal F}$ one can define a $G$-invariant
connection in a natural way, whose curvature is denoted by ${\bf
K}$. The $G$-action on ${\mathcal F}$ has a moment map
$\mu:{\mathcal F}\to \frak{gl}(H)\otimes {\frak g}^*$, relative to
the $2$-form ${\bf K}$; that is, $D\langle\mu,\,A\rangle=-{\bf
K}(Y_A,\,.),$ where $D$ is the covariant derivative. Moreover
$\langle\mu,\,A\rangle=:{\bf h}_A$ induces a map $h_A$ from
${\mathcal O}$ to $\frak{gl}(H)$.

 We will denote by ${\mathcal V}$  the vector bundle
with fibre $H$ associated to ${\mathcal F}$.
 We write  ${\mathcal B}_2$ for the space of smooth sections $\sigma$ of
 ${\mathcal V}$ which can be identified with   the elements of ${\mathcal
 B}_1$, and
  %with compact support and polarized with respect to ${\frak
  %u}$. This space of sections can be identified in a natural way
 % with ${\mathcal B}_1$.
  we put ${\mathcal B}_3$ for the space of maps $f:{\mathcal
 F}\to H$ associated with the sections of ${\mathcal V}$ that belong  to ${\mathcal B}_2$.

  Given $A\in {\frak g}$ we will denote by $X_A$ the vector field
  on ${\mathcal O}$ defined by the coadjoint action of $G$. On the space ${\mathcal B}_2$
  %sections of ${\mathcal V}$,
   we consider the following operator
$${\mathcal P}_A:=-D_{X_A}+h_A.$$
In Section \ref{S:derivative} we prove the following theorem, that
gives the representation $\pi'$ on the spaces ${\mathcal B}_i$,
for $i=2,3$.

\begin{Thm}\label{theorempiprime} The  representations of
${\frak g}$
$$A\in{\frak g}\mapsto {\mathcal P}_A\in{\rm End}({\mathcal B}_2)$$
and
$$A\in{\frak g}\mapsto  -Y_A\in{\rm End}({\mathcal B}_3)$$
are equivalent to $\pi'$, the differential representation of
$\pi$.
 \end{Thm}

Theorem \ref{theorempiprime} gives the representation $\pi'$ on
geometric objects. To determine a geometric description of
$\pi(g_1)$ we will ``integrate" $\pi'$ along a curve in $G$ with
final point at $g_1$. To abbreviate,   the smooth curves    in $G$
with initial point at $e$ will be called {\em paths} in $G$.
 Let $\{g_t\}$ be a path in $G$
 %initial point at $e$ and
with $g_1$ in $Z(G)$, the center of $G$. This curve determines its
velocity curve; that is, the family $\{A_t\}\subset{\frak g}$
given by the relation
 $\Dot g_tg_t^{-1}=A_t.$
The corresponding time-dependent vector field $Y_{A_t}$ defines a
Hamiltonian  flow ${\bf F}_t$ on ${\mathcal F}$. We will prove
that the time-$1$ map of this flow is precisely the multiplication
by $\kappa$ in ${\mathcal F}$ (see item (d) of Theorem
\ref{kappatransport}).

Given $s\in{\mathcal B}_1$, we define a family of maps $\{s_t:
G\to H\}_{t}$  by the equations
 \begin{equation}\label{trnsportG}
  \frac{d\,s_t}{dt}=-R_{A_t}s_t,\;\;\; s_0=s.
 \end{equation}

Given a section $\sigma\in{\mathcal B}_2$ we  consider the family
of sections $\sigma_t$ of ${\mathcal V}$ determined by the
following equations
\begin{equation}\label{fracdsigmat}
\frac{d\,\sigma_t}{dt}={\mathcal P}_{A_t}\sigma_t,\;\;\;
\sigma_0=\sigma.
 \end{equation}

Similarly, given $f\in{\mathcal B}_3$, let $\{f_t\}$ be the set of
maps $f_t:{\mathcal F}\to H$ such that
\begin{equation}\label{definft}
\frac{d\,{f}_t(p)}{dt}=-Y_{A_t}(p)(f_t),\;\;\;\;
  f_0=  f.
  \end{equation}

  The following theorem
 %that gives geometric interpretations of
 relates
 the constant $\kappa$ with
  the solutions of the ``evolution" equations (\ref{trnsportG}),
   (\ref{fracdsigmat}), (\ref{definft}) and with the time-$1$ map
   ${\bf F}_1$.

 \begin{Thm}\label{kappatransport}
Let ${g_t}$ be an arbitrary path on $G$ with $g_1\in Z(G)$, and
$A_t$ the corresponding velocity curve. If $\kappa$ is the
constant given by
%If $g_1\in Z(G)$, let $\kappa$ be the constant
  %given by
%Proposition \ref{g1Z(G)}
(\ref{kappadefinition}), the following statements hold

 (a) If $s_t$ is the solution of (\ref{trnsportG}), then
 $s_1=\kappa s$.

 (b) If $\sigma_t$ is the solution to (\ref{fracdsigmat}), then
$\sigma_1=\kappa\sigma.$

 (c)  If $f_t$   is the solution to (\ref{definft}), then
$f_1=\kappa f.$

 (d) ${\bf F}_1$ is the multiplication by $\kappa$; that is, ${\bf
F}_1[g,\,\alpha]=\kappa[g,\,\alpha]$.
\end{Thm}

For the above   path $\{g_t\}$ with endpoint at $g_1\in Z(G)$, we
denote by $\psi_t$ the closed isotopy on ${\mathcal O}$ determined
by the time-dependent vector field $X_{A_t}$, that is,

\begin{equation}\label{isotopy}
\frac{d\,\psi_t}{dt}=X_{A_t}\circ\psi_t,\;\;\; \psi_0={\rm id}.
\end{equation}

When ${\rm dim}\,H=1$ the curvature ${\bf K}$ projects a $2$-form
$-\omega$ on the orbit ${\mathcal O}$, and  ${\bf F}_1$ is the
  action integral (\cite{We}, \cite{Mc-S}) around $\psi_t$ (see Section \ref{Sect.Physical}). This
fact is the statement of the following theorem, which will be
proved in Section \ref{S:Isotopies}
\begin{Thm}\label{integralaccion} If ${\rm dim}\,H=1$
 $$\kappa={\rm
exp}\Big(\int_S\omega+\int_0^1 h_{A_t}(q_t)dt \Big),$$ where $q$
is an arbitrary point of ${\mathcal O }$ and $S$ any $2$-chain in
${\mathcal O}$ whose boundary is the curve $\{q_t:=\psi_t(q)\}_t$.
 \end{Thm}

Let $K$ be a maximal compact subgroup of $G$, and let $\tilde g$
an element of $K$, such that its conjugacy class meets $L_0$, the
connected component of the identity of $L$; that is,
 there exists $a\in G$ such that $a^{-1}\tilde g a\in L_0$. Let $g_t$
 be a path in $G$ with $g_1=\tilde g$ and  $a^{-1} g_t a\in L$,
  $\,A_t$  the corresponding velocity path and
 $h_{A_t}$ the map on ${\mathcal O}$ induced by $\langle
 \mu,\,A_t\rangle$. Let $\tau$ be an irreducible representation of
 $K$ which occurs in the representation $\pi_{|K}$ with
 nonzero multiplicity.
%In the case that ${\rm dim}\,H=1$,
 We will prove the
 following theorem, that gives an expression for the value of
 character $\chi_{\tau}$ at $\tilde g$ in terms of the
 ``Hamiltonian" functions $h_{A_t}$.

 \begin{Thm}\label{theoremcharacter}
Let $\tilde g$ be an element of $K$, such that there exists $a\in
G$ with $a^{-1}\tilde g a$ in $L_0$.  Let $C_t$ be the velocity
curve of a path in $L$ with endpoint at $a^{-1}\tilde g a$. If
${\rm dim }\,\tau=m$ and ${\rm dim}\,H=1$, then
 $$\chi_{\tau}(\tilde g)=m\,{\rm
exp}\Big(\int_0^1 h_{A_t}(x_0) dt\Big),$$
 where $A_t={\rm Ad}_a
C_t$ and $x_0=a\cdot\eta\in{\mathcal O}.$
\end{Thm}

 An important invariant of the representation $\pi'$  is its
 infinitesimal character $\chi$ defined on ${\mathcal Z}({\frak g}_{\mathbb C})$, the
 center of the universal enveloping algebra $U({\frak g}_{\mathbb C})$ of
 ${\frak g}_{\mathbb C}={\frak g}\otimes {\mathbb C}\,$ \cite{Kn}.
 If  ${\frak h}$ is a Cartan subalgebra of ${\frak g}_{\mathbb C}$
 contained in ${\frak l}_{\mathbb C}$ and $Z\in {\mathcal Z}({\frak g}_{\mathbb
 C})$, we will denote by  $\Hat Z$  the projection of $Z$
 into
 $U({\frak h})$. Let
 $Y_1,\dots, Y_r$ be a basis of ${\frak h}$.
For each $Y_i$ we define the map ${\bf h}_i:{\mathcal
 F}\to \frak{gl}(H)$, by ${\bf h}_i=\langle\mu,\,Y_i\rangle$; that
 is, ${\bf h}_i$ is the Hamiltonian function associated with $Y_i$.
 If $q(Y_1,\dots,Y_r)$ is a   polynomial
\begin{equation}\label{Polynomial}
 q(Y_1,\dots,Y_r)=a+\sum_k
 a_kY_k+\sum_{i,j}a_{ij}Y_iY_j+\sum_{k,i,j}a_{kij}Y_kY_iY_j+\dots\;\;
 \hbox{(finite sum)}
 \end{equation}
 such that $\Hat Z=q(Y_1,\dots,Y_r)\in  U({\frak
 h})$, a geometric interpretation of $\chi(Z)$ is given in
 the following theorem

\begin{Thm}\label{7rule}
If   $Z$ is an element of ${\mathcal Z}({\frak g}_{\mathbb C})$
such that $\Hat Z$ is defined by the   polynomial $q(Y_1,\dots,
Y_r)$. Then
$$q({\bf h}_1,\dots,{\bf h}_r):{\mathcal F}\to  \frak{gl}(H)$$
is a constant map on the fiber over $\eta$, and its value on this
fiber is $\chi(Z) {\rm Id}$.
% the infinitesimal character of $\pi'$ on $Z$.
\end{Thm}

%%%%%%%%%%%%%%%%%%%%%%%%%%%%%%%%%%%%%%%%%%%%%%%%%%%%%%%%%%%%%%%%%%%%%%%%%%%%%%%%%%

% This result
 % can be considered as a new version of the rule 7 of the Orbit Method proposed by
% Kirillov in \cite{Kir}, page xix.

%%%%%%%%%%%%%%%%%%%%%%%%%%%%%%%%%%%%%%%%%%%%%%%%%%%%%%%%%%%%%%%%%%%%%%%%%%%%%%%%%%%%%%%%%%%%%%%%%%%%%%%%%%%%%%%%%%%%%%%%%%%%%%%%%

\smallskip

 This article  is organized as follows. In Section \ref{S:definitions}
 we introduce the definitions and notations which will be used. Following Vogan  (\cite{Vo1}, \cite{Vo2})
 we define the representation $\pi$ associated to the orbit ${\mathcal O}$ of a hyperbolic element.

 In Section \ref{S:derivative} we describe the differential representation $\pi'$ on the spaces
 ${\mathcal B}_i$, for $i=2,3$, proving Theorem
 \ref{theorempiprime}.

In Section \ref{S:Isotopies} we give geometric interpretations of
 $\pi(g_1)$, for $g_1\in Z(G)$. We will prove   Theorem
 \ref{kappatransport} and Theorem \ref{integralaccion}.  In
 Subsection \ref{Subsection:the character} we prove Theorem \ref{theoremcharacter} about
  the character of $\pi$. Subsection
 \ref{Subsection:infinitesimal}
 concerns with the geometric
interpretation of the infinitesimal character of $\pi'$; in this
subsection we will prove Theorem \ref{7rule}.

Section \ref{Sect.Physical} provides an interpretation of Theorem
\ref{integralaccion} in terms of physical concepts. We will show
  that the invariant $\kappa$ can be considered as the exponential of the action
integral around the closed curve $\psi_t$,
 %in the Hamiltonian group of ${\mathcal O}$
   and also as the Berry phase of a loop of Lagrangian
   submanifolds of ${\mathcal O}$. In a worked example we will
   consider a hyperbolic orbit of the restricted Lorentz group
   $SO^+(1,3)$. Using Theorem \ref{7rule} we will calculate the value of the corresponding
   infinitesimal character on the Casimir element $C$, and we will
   interpret this value in terms of the ``quantum" operator that Geometric Quantization
   associates with   $C$.

\smallskip

 %%%%%%%%%%%%%%%%%%%%%%%%%%%%%%%%%%%%%%%%%%%%%%%%%%%%%%%%%%%%%%%%%%%%%%%%%%%%%%%%%%%%%%%%%%
 %%%%%%%%%%%%%%%%%%%%%%%%%%%%%%%%%%%%%%%%%%%%%%%%%%%%%%%%%%%%%%%%%%%%%%%%%%%%%%%%%%%%%%%%%%%

\section {Definitions and notations.} \label{S:definitions}

 %As we said,
 Here we review the construction of the representation
associated to the coadjoint orbit of a hyperbolic element (see
\cite{Vo1}, \cite{Vo2} for details).

% Let $G$ be a reductive group. As definition of reductive group   we adopt the one
% given by Vogan in \cite{Vo2}. If $g\in G$ and $A\in{\frak g}$, we
% put $g\cdot A={\rm Ad}_g(A)$, and if $\xi\in{\frak g}^*$ we write
% $g\cdot\xi$ for ${\rm Ad}^*_g(\xi)$.

By $G$ we denote a reductive group. As definition of reductive
group we adopt the one given by Vogan in \cite{Vo2}. We recall
this definition. A {\em linear} group is reductive if it has
finitely many connected components and is preserved by the Cartan
involution. A reductive group $G$ is a Lie group endowed with a
homomorphism from $G$ onto $G_1$ of finite kernel, $G_1$ being a
linear reductive group. In particular, $GL(n,{\mathbb R})$,
$SL(n,\,{\mathbb R})$, $SO(p,q)$, $O(p,q$), $Sp(2n)$  and all
compact Lie groups are reductive groups.
 If $g\in G$ and $A\in{\frak g}$, we
 put $g\cdot A={\rm Ad}_g(A)$, and if $\xi\in{\frak g}^*$ we write
 $g\cdot\xi$ for ${\rm Ad}^*_g(\xi)$.

Let ${\frak g}={\frak k}\oplus{\frak p}$ be a Cartan decomposition
of ${\frak g}$, with ${\frak k}$ the Lie algebra of a maximal
compact subgroup $K$ of $G$.

Given $\eta\in{\frak g}^*$, its stabilizer  for the coadjoint
adjoint of $G$ will be denoted by $L$.  On the other hand, $\eta$
determine an element
 $X_0\in{\frak g}$ by the
equality
  \begin{equation}\label{X0}
  \eta(Y)={\rm Re \; Tr}(X_0Y),\;\;\hbox{for all}\;
 Y\in{\frak g},
  \end{equation}
   where ${\frak g}$ is identified with a Lie algebra
 of matrices.

Let us assume that $\eta$ is a {\em hyperbolic} element of ${\frak
g}^*$. We can suppose that $X_0\in{\frak p}$, after replacing
$\eta$ by an ${\rm Ad}^*(G)$-equivariant element. As ${\rm
ad}(X_0)$ is a diagonalizable endomorphism of ${\frak g}$ with
real eigenvalues
$${\frak g}=\bigoplus_{r\in{\mathbb R}}{\frak g}_r,\;\;\;
{\frak g}_r=\{Y\in{\frak g}\,|\, {\rm ad}(X_0)(Y)=rY   \},$$
 and ${\frak g}_0$ is the Lie algebra ${\frak l}$ of the
 subgroup $L$. Moreover the adjoint action of $l\in L$
 preserves each ${\frak g}_r$.

 We put
  \begin{equation}\label{fraku}
{\frak u}=\bigoplus_{r>0}{\frak g}_r,\;\;\; {\frak
u}^-=\bigoplus_{r<0}{\frak
  g}_r,\;\;\; U={\rm exp}\,{\frak u}.
   \end{equation}

  Then $U$ is a simply connected nilpotent subgroup of $G$ and $L$
  normalizes $U$. So $Q:=LU$ is a Levi decomposition of the
  subgroup $Q$.

  We define the following positive character on $Q$
  \begin{equation}\label{Delta}
   \Delta(q)=|\,{\rm det}\,({\rm Ad}(q)_{|{\frak u}})|^{1/2}.
   \end{equation}
   The derivative of $\Delta$ will be denoted by $\delta$. Since
  $[{\frak g}_r,\,{\frak g}_s]\subset {\frak g}_{r+s}$,
  $\delta(B)=0$ for all $B\in {\frak u}$. On the other hand
  $\Delta(l\,l_1\,l^{-1})=\Delta(l_1)$, for all $l,\,l_1\in{L}$,
  so $\delta(l\cdot A)=\delta(A)$, for $l\in L$ and $A\in{\frak
  l}$.
 We extend $\delta$ to a linear map on ${\frak g}$ by setting
  $\delta_{|{\frak u}^-}=0$. As the action of $ L$ preserves
  the ${\frak g}_r$, then
 % Hence $\delta[A,\,Y]=0$ for $A\in{\frak
  %l}, \,Y\in{\frak g}$. Consequently
  \begin{equation}\label{invdelta}
  \delta(l\cdot Y)=\delta(Y),\;\; \hbox{for all}\; l\in L \;\;\hbox{and}\; Y\in{\frak g}.
  \end{equation}

Let $\Lambda$ be an {\em integral datum} at $\eta$ \cite{Vo1};
that is, $\Lambda$ is an irreducible unitary representation of $L$
in a Hilbert space $H$, such that
 \begin{equation}\label{datum}
 \Lambda({\rm exp}\,A)=e^{i\eta(A)}{\rm Id},\;\;\hbox{for all}\; A\in{\frak l}.
  \end{equation}
 So $\Lambda(l\,{\rm exp}( A)\,l^{-1})=\Lambda({\rm exp}\, A)$.
We extend $\Lambda$ to $Q$ by $\Lambda(lu)=\Lambda(l)$, and write
$\lambda$ for the  derivative of $\Lambda$.   In turn, $\lambda$
can be extended to a linear map on ${\frak g}$ by putting
$\lambda_{|{\frak u}^-}=0$. As in the preceding case
 \begin{equation}\label{invlambda}
  \lambda(l\cdot Y)=\lambda(Y),\;\; \hbox{for all}\; l\in L \;\;\hbox{and}\; Y\in{\frak g}.
  \end{equation}

${\rm GL}(H)$ will denote the group of continuous linear
automorphisms of $H$. We put $\Phi:Q\to {\rm GL}(H)$ for the
representation of $Q$ tensor product $\Lambda\otimes \Delta$. By
$\phi$ we denote the  linear map
 \begin{equation}\label{definphi}
  \phi:A\in{\frak g}\mapsto
\lambda(A)+\delta(A){\rm Id}\in \frak {gl}(H).
 \end{equation}
 From (\ref{invdelta}) and (\ref{invlambda})  it follows
 \begin{equation}\label{invphi}
 l\cdot\phi=\phi.
 \end{equation}

  By $\pi$ we denote the representation ${\rm Ind}_Q^G(\Phi)$; that is,    the irreducible unitary
  representation of $G$ induced by $\Phi$. The space of   $C^{\infty}$
  vectors of $\pi$ is the space of smooth   functions $s:G\to H$,
  with compact support modulo $L$,
  such that
  \begin{equation}\label{equivariant}
s(g\,l)=\Phi(l^{-1})s(g),\;\;\hbox{for all}\,\, g\in G,\; l\in L;
 \;\;\;\;\hbox{and}\;\; L_A s=0,\;\; \hbox{for all}\,\, A\in
{\frak u},
  \end{equation}
where $L_A$ is the left invariant vector field on $G$ defined by
$A$. The representation $\pi$ on this space is given by
\begin{equation}\label{equivariantbis}
 \pi(g)(s)=s\circ L_{g^{-1}},
 \end{equation}
  $L_{g^{-1}}$ being the left
multiplication in $G$ by $g^{-1}$. Therefore the differential
representation   $\pi'$ of $\pi$ on the smooth function $s$, with
compact support modulo $L$, which satisfy (\ref{equivariant}) is
given by
\begin{equation}\label{piprime}
\pi'(C)(s)=-R_C(s),
\end{equation}
$R_C$ being the right invariant vector filed on $G$ defined by
$C\in\frak{g}$.

%%%%%%%%%%%%%%%%%%%%%%%%%%%%%%%%%%%%%%%%%%%%%%%%%%%%%%%%%%%%%%%%%%%%%%%%%%%%%%%%%%%%%%%%%%%%%%%%%%%%%%%%%%%
%%%%%%%%%%%%%%%%%%%%%%%%%%%%%%%%%%%%%%%%%%%%%%%%%%%%%%%%%%%%%%%%%%%%%%%%%%%%%%%%%%%%%%%%%%%%%%%%%%%%%%%%%%%%%%

%\section{Principal fibre bundle}

\section{The differential representation}\label{S:derivative}

Henceforth
 % and until  Section \ref{SectionOrbit},
% Orbit of an elliptic element In Section \ref{S:infinitesimal}
 $G$ will be a reductive group, $\eta$ a hyperbolic element of
${\frak g}^*$ and $\Lambda$ an integral datum at $\eta$. The
coadjoint orbit of $\eta$ will be denoted by ${\mathcal O}$. If
$B\in\frak{g}$, $\,X_B(g\cdot\eta)$ will be the tangent vector to
${\mathcal O}$ at $g\cdot \eta$ defined by the curve
 $\{{\rm exp}(tB)\cdot(g\cdot\eta) \}_t$.
 %With $B^*$ we denote
 %the vector field on ${\mathcal O}$ defined by the curves
%$\{ g\cdot{\rm exp}(tB)\cdot\eta) \}_t$.

By (\ref{invphi}) the map
 \begin{equation}\label{definh}
  h_B:g\in
G\mapsto \phi(g^{-1}\cdot B)\in\frak{gl}(H)
 \end{equation}
 induces a mapping on
$\mathcal{O}$, that will be also denoted by $h_B$. For any
$A,\,B\in{\frak g}$
\begin{equation}\label{XA(hA)}
X_A(h_B)=-h_{[A,\,B]}.
\end{equation}

Next we define the following $GL(H)$-principal bundle over
$\mathcal{O}\simeq G/L$.
$$\mathcal{F}=G\times_L GL(H)=\{ (g,\,\alpha)\,|\, g\in G,\, \alpha\in
GL(H)\}/\sim,$$
 where $(g,\,\alpha)\sim (gl,\,\Phi(l^{-1})\alpha)$, with $l\in
 L$.

 On $\mathcal{F}$ there is a natural left $G$-action ${\mathcal L}$,
  so   each $B\in\frak{g}$ determines  a vector field $Y_B$ on $\mathcal{F}$ and
  $$({\mathcal L}_g)_*(Y_B)=Y_{g\cdot B}.$$

On the other hand,  each $y\in\frak{gl}(H)$ determines a vertical
vector field $W_y$ on  $\mathcal{F}$ by means  of right
$GL(H)$-action ${\mathcal R}$.

 Given $A\in{\frak l}$,  the trivial curve $\{[g\cdot e^{tA},\, \Phi(e^{-tA})\alpha ]
 \}_t$ in ${\mathcal F}$ defines the vector $Y_{g\cdot A}[g,\,\alpha]-W_y[g,\,\alpha]$, with $y={\rm
 Ad}\,\alpha^{-1}(\phi(A))$, where ${\rm Ad}$ denotes the adjoint
action of $GL(H)$ on $\frak{gl}(H)$. So the tangent space to
$\mathcal{F}$ at $[g,\,\alpha]$ is
 \begin{equation}\label{tangentF}
T_{[g,\alpha]}\mathcal{F} =\frac{ \{ Y_B[g,\,\alpha] \;|\;
B\in{\frak g} \}\oplus\{W_y[g,\,\alpha] \;|\;y\in\frak{gl}(H) \} }
  { \{ Y_{g\cdot A}[g,\,\alpha]-W_{{\rm Ad}\,\alpha^{-1}(\phi(A))  }[g,\,\alpha]\;|\; A\in{\frak l}  \}    }
 \end{equation}

On ${\mathcal F}$ we define the following $\frak{gl}(H)$-valued
$1$-form
\begin{equation}\label{definOmega}
\Omega(Y_B[g,\,\alpha]+W_y[g,\,\alpha])={\rm
Ad}\,\alpha^{-1}(\phi(g^{-1}\cdot B))+y.
\end{equation}

 $\Omega$ is, in fact, well-defined on the quotient
(\ref{tangentF}), and  it is easy to check that
$${\mathcal L}_g^*\Omega=\Omega, \;\;\hbox{and}\;\;
 {\mathcal R}_{\alpha}^*\Omega={\rm Ad}_{\alpha^{-1}}\circ\Omega,$$
 for $g\in G$ and $\alpha\in GL(H)$.
 That is, one has the proposition
 \begin{Prop}\label{Proconnection}
 The $1$-form $\Omega$ defined in (\ref{definOmega}) is a $G$-invariant connection on the $GL(H)$-principal bundle
 ${\mathcal F}$.
 \end{Prop}

 We can lift $h_A$ to a a well-defined function ${\bf h}_A:{\mathcal F}\to \frak
{gl}(H)$ by setting
 \begin{equation}\label{lifhA}
 {\bf h}_A[g,\,\alpha]={\rm
Ad}_{\alpha^{-1}}h_A(g).
 \end{equation}
  Then
\begin{equation}\label{defiOmega}
\Omega(Y_B+W_y)={\bf h}_B+y
 \end{equation}

\begin{Lem}\label{Lembfh} Given $A\in{\frak g}$ and
$y\in\frak{gl}(H)$, then
 $$Y_A({\bf h}_B)=-{\bf h}_{[A,\,B]}
 \;\;\,\hbox{and}\;\;\, W_y[g,\,\alpha]({\bf h}_B)=-\big[y,\,{\bf
h}_B([g,\,\alpha])\big]_{\frak{gl}},$$
 where $[\; ,\;]_{\frak
{gl}}$ is the bracket in the Lie algebra $\frak
 {gl}(H)$.
 \end{Lem}
{\it Proof.} From (\ref{XA(hA)}) it follows
$$Y_A[g,\,\alpha]({\bf h}_B)=
 \frac{d}{dt}\bigg|_{t=0}{\bf h}_B[{\rm exp}(tA)\cdot g,\,\alpha]=
{\rm Ad}_{\alpha^{-1}}\big(X_A(g)(h_B)\big)=-{\bf
h}_{[A,\,B]}[g,\,\alpha].$$

 The second formula can be directly deduced from
 $${\bf h}_B[g,\,\alpha\,{\rm exp}(ty)]={\rm Ad}\,_{{\rm exp}(-ty)}{\bf
 h}_B[g,\,\alpha].$$

 \qed

 Next we will calculate the value of the curvature ${\bf K}$ of
 the
 connection  $\Omega$ on the pair $(Y_B,\,Y_C)$ of vector fields.
\begin{Prop}\label{ProbfK}
The curvature ${\bf K}$ of the connection $\Omega$ satisfies
\begin{equation}\label{Curvature}
 {\bf K} (Y_B,\,Y_C)=
  -{\bf h }_{[B,\,C]} +\big[{\bf h}_B ,\,{\bf h}_C \big]_{\frak
  {gl}},
 \end{equation}
for all $B,C\in{\frak g}$.
\end{Prop}

{\it Proof.}
  By the
 structure equation
 $${\bf K}(Y_B,\,Y_C)=d\,\Omega(Y_B,\,Y_C)+\big[\Omega(Y_B),\,\Omega(Y_C)  \big]_{\frak{gl}}.$$
 From Lemma \ref{Lembfh}  and (\ref{defiOmega}) it follows
%\begin{align}\label{aligniacion}
 $$Y_B[g,\,\alpha](\Omega(Y_C))= {\bf
h}_{[C,\,B]}[g,\,\alpha].$$
 Similarly $Y_C(\Omega(Y_B))={\bf h}_{[B,\,C]}.$
  % {\rm Ad}\,\alpha^{-1}(\phi(g^{-1}\cdot{\rm exp}(-tB)\cdot
 %C)) \big)= \\ \notag
 %=&{\rm Ad}\,\alpha^{-1}(\phi(g^{-1}\cdot[C,\,B]))={\rm
 %Ad}\,\alpha^{-1}(h_{[C,B]}(g)).
 % \end{align}
 %\begin{align}\label{aligniacion}
 %Y_B[g,\,\alpha](\Omega(Y_C))=&
 %\frac{d}{dt}\bigg|_{t=0}\big( {\rm Ad}\,\alpha^{-1}(\phi(g^{-1}\cdot{\rm exp}(-tB)\cdot
 %C)) \big)= \\ \notag
% =&{\rm Ad}\,\alpha^{-1}(\phi(g^{-1}\cdot[C,\,B]))={\rm
% Ad}\,\alpha^{-1}(h_{[C,B]}(g)).
 %\end{align}
Hence $d\,\Omega(Y_B,\,Y_C) =- {\bf h}_{[B,C]}$ and
(\ref{Curvature}) follows.

 \qed

By $D$ we denote the covariant derivative determined by $\Omega$.
Since the horizontal component of $Y_A[g,\,\alpha]$ is
$Y_A[g,\,\alpha]+W_y[g,\,\alpha]$, with $y=-{\bf
h}_A[g,\,\alpha]$, by Lemma \ref{Lembfh}
$$
 \notag D\,{\bf h}_B(Y_A[g,\,\alpha]) = -{\bf h}_{[A,\,B]}+\big[{\bf
 h}_A, \,{\bf h}_B
 \big]_{\frak{gl}}$$

 It follows from (\ref{Curvature}) that
$$D\,{\bf h}_B(Y_A)=-{\bf K}(Y_B,\,Y_A),\;\;\hbox{for all}\; A\in{\frak
g}.$$
 Thus we have
 \begin{equation}\label{DHathB}
 D\,{\bf  h}_B=-\iota_{Y_B}{\bf K},\;\;\hbox{for all}\; B\in{\frak
g}.
 \end{equation}

 Equation (\ref{DHathB}) can be interpreted
saying that the $\frak{gl}(H)$-valued function ${\bf  h}_B$
  is a ``Hamiltonian" for the vector field $Y_B$, with respect to
the covariantly closed
   $\frak{gl}(H)$-valued $2$-form ${\bf K}$. That is, we define
   $$\mu:{\mathcal F}\to\frak{gl}(H)\otimes\frak{g}^*$$
     by $\langle\mu,\,A\rangle={\bf h}_A$, with $A\in{\frak g}$. This map is
    $G$-equivariant, that is, $\langle
    \mu(gp),\,A\rangle=\langle\mu(p),\,g^{-1}\cdot A\rangle$, for
    all $p\in{\mathcal F}$ and all $g\in G$ and
    $D\langle\mu,\,A\rangle=-\iota_{Y_B}{\bf K}$. We call $\mu$ the
    moment map for the $G$-action on $({\mathcal F},\,{\bf K})$.

\smallskip

We denote by ${\mathcal V}$ the vector bundle on ${\mathcal O}$
with fibre $H$
  $$G\times_L H=\{\langle g,\,v\rangle \,|\, g\in
G,\,v\in H\},$$
 with
$\langle g,\,v\rangle=\langle gl,\,\Phi(l^{-1})v\rangle.$ The
vector bundle ${\mathcal V}$ can also be considered as associated
to ${\mathcal F}$ by the natural representation of $GL(H)$. That
is, ${\mathcal V}$ is
 $$\big\{ \{p,\,v\}\,|\,p\in{\mathcal F},\, v\in H \big\},$$
 where $\{p,\,v\}=\{p\,\beta,\,\beta^{-1}v\}$ for all
$\beta\in GL(H)$. The correspondence $\langle g,\,
v\rangle\mapsto\{[g,\,{\rm Id}],\,v\}$ gives the isomorphism
between those definitions of ${\mathcal V}$.

Now  we consider the following three vector spaces of {\em smooth}
maps
%$${\mathcal B}_1=\{s:G\to H\,|\, s\;\hbox{smooth map satisfying}\;
%(\ref{equivariant})\}$$
%$${\mathcal B}_2=\{ \tau\,|\, \tau\;\hbox{smooth section of}\; {\mathcal V}\;\hbox{satisfying}\;
%D_{A^*}\tau=0,\, \forall A\in{\frak u}\}$$
%$${\mathcal B}_3=\{f:{\mathcal F}\to H \;\hbox{smooth} \,|\,
%f(p\,\beta)=\beta^{-1}f(p)\;\forall p\in {\mathcal F}, \beta\in
%GL(H);  Z_{A}f=0,\, \forall A\in{\frak u}\}.$$
 $${\mathcal B}'_1=\{s:G\to H\,|\,
 s(gl)=\Phi(l^{-1})s(g),\,\forall g\in G,\,\forall l\in L,\,{\rm supp}(s)\,\, \text{compact modulo}\, L  \}$$
$${\mathcal B}'_2=\{ \tau\,|\, \tau\;\hbox{section of}\;\, {\mathcal V},\ \hbox{supp}\,(\tau)\; \hbox{compact}\}$$
 $${\mathcal B}'_3=\{f:{\mathcal F}\to H  \,|\,
f(p\,\beta)=\beta^{-1}f(p),\,\forall p\in {\mathcal F},
\forall\beta\in GL(H),\,{\rm pr}({\rm supp}(f))\; \hbox{compact}
\},$$ where ${\rm pr}$ is the projection ${\rm pr}:{\mathcal F}\to
G/L$.

Given $s\in {\mathcal B}'_1$, determines a section $\sigma\in
{\mathcal B}'_2$ by the relation
\begin{equation}\label{sigma-s}
\sigma(gL)=\langle g, \, s(g)\rangle.
 \end{equation}
   Moreover $s$ defines
$\sigma^{\sharp}\in{\mathcal B}'_3$ by
 \begin{equation}\label{sigmasaharps}
\sigma^{\sharp}[g,\,\alpha]=\alpha^{-1}s(g).
 \end{equation}
  With the above notations
\begin{equation}\label{sigm-sigmasharp}
\sigma(x)=\{p,\,\sigma^{\sharp}(p)\},
 \end{equation}
 for any
$p\in{\mathcal F}$ in the fibre of $x\in{\mathcal O}.$

It is well-known that the correspondences $s\mapsto\sigma$ and
$s\mapsto\sigma^{\sharp}$ allow us to identify the ${\mathcal
B}'_i$'s. We denote also by $D$ the covariant derivative on
sections of ${\mathcal V}$ defined by the connection $\Omega$. It
is a known fact that the map of ${\mathcal B}'_3$  associated with
the section $D_{X_A}\sigma$ of ${\mathcal V}$ is
$X^{\sharp}_A(\sigma^{\sharp})$, where $X^{\sharp}_A$ is the
horizontal lifting of the vector field $X_A$.

\smallskip

We set
 \begin{equation}\label{B1}
{\mathcal B}_1=\{s\in{\mathcal B}'_1\,|\,L_As=0,\;\forall
A\in{\frak u}\}.
 \end{equation}
 Since $Q$ is
connected the condition defining ${\mathcal B}_1$ is equivalent to
$s(gq)=\Phi(q^{-1})s(g)$ for all $q\in Q$ and all $g\in G$. In
order to interpret ${\mathcal B}_1$ in terms of sections of
${\mathcal V}$ and equivariant functions on ${\mathcal F}$ we need
an additional fibre bundle.

${\mathcal F}_Q$ is the $GL(H)$-principal fibre bundle over $G/Q$
defined by $\Phi$; that is,
 $${\mathcal F}_Q:=G\times_Q GL(H).$$
 % The class in ${\mathcal
%F}_Q$ of $(g,\,\alpha)$ will be denoted $[g,\,\alpha]_Q$.
 One has
a natural fibre map $\Xi:{\mathcal F}\to{\mathcal F}_Q$ over the
canonical projection
\begin{equation}\label{Projection}
 G/L\to G/Q.
 \end{equation} We put
  \begin{equation}\label{B3}
{\mathcal B}_3= \{f\in{\mathcal B}'_3\,|\, f\; \text{ factors
through }\; \Xi \}.
 \end{equation}

Analogously we define ${\mathcal V}_Q:=G\times_{Q}H$, and the
natural fibre map $\Xi:{\mathcal V}\to{\mathcal V}_Q$ will be also
denoted by $\Xi$. We set
  \begin{equation}\label{B2}
 {\mathcal B}_2:=\{\tau\in{\mathcal
B}'_2\,|\, \Xi\circ\sigma\; \text{is constant along the fibers of
}\; (\ref{Projection}) \}.
\end{equation}

{}From the above definitions is easy to prove the following
proposition
\begin{Prop}\label{PropBi}
  The correspondences $s\mapsto \sigma$ and
  $s\mapsto\sigma^{\sharp}$ define bijective maps between the ${\mathcal
  B}_i$'s.
\end{Prop}

{\bf Proof of Theorem \ref{theorempiprime}.} From
(\ref{sigmasaharps}) one deduces
\begin{equation}\label{YAWy}
Y_A[g,1](\sigma^{\sharp})=R_A(g)(s),\;\;\hbox{and}\;\;
 W_y[g,1](\sigma^{\sharp})=-y\,s(g).
\end{equation}
On the other hand, by (\ref{defiOmega})
\begin{equation}\label{Xsharp}
X^{\sharp}_B[g,\,\alpha]=Y_B[g,\,\alpha]+W_y[g,\,\alpha],
\end{equation}
with $y+ {\bf h}_B[g,\,\alpha]=0$.

 From (\ref{YAWy}) together with (\ref{Xsharp})
 %,(\ref{sigm-sigmasharp}) and (\ref{sigma-s})
 and (\ref{sigmasaharps}) it follows that the
$\Phi$-equivariant function on $G$ associated with $D_{X_A}\sigma$
is $R_A(s)+h_A\,s$. So the section $-D_{X_A}\sigma+h_A\sigma$ of
${\mathcal V}$ has as associated $\Phi$-equivariant function to
$-R_A(s)$.
  Therefore if
we put
  $${\mathcal P}_A(\sigma):=-D_{X_A}\sigma+h_A\sigma,$$ then the
  family $\{{\mathcal P}_A\}$ of endomorphisms  is a representation of ${\frak g}$ on ${\mathcal
  B}_2$ equivalent to $\pi'$ defined in (\ref{piprime}).

  From (\ref{sigmasaharps}) it follows
$$Y_A[g,\,\alpha](\sigma^{\sharp})=\alpha^{-1}R_A(g)(s).$$
That is, $Y_A\,\sigma^{\sharp}$ is the function of ${\mathcal
B}_3$ associated to $R_A\,s\in{\mathcal B}_1$. Hence the algebra
representation $\pi'$ defined in (\ref{piprime}) is equivalent to
the representation of ${\frak g}$ on the space ${\mathcal B}_3$
given by the
 %vector fields
 operators $\{-Y_A\}$.

\qed

%%%%%%%%%%%%%%%%%%%%%%%%%%%%%%%%%%%%%%%%%%%%%%%%%%%%%%%%%%%%%%%%%%%%%%%%%%%%%%%%%%%%%%%%%%%%%%%%%%%%%%%%%%%%%%%%%%%%%%%%
%%%%%%%%%%%%%%%%%%%%%%%%%%%%%%%%%%%%%%%%%%%%%%%%%%%%%%%%%%%%%%%%%%%%%%%%%%%%%%%%%%%%%%%%%%%%%%%%%%%%%%%%%%%%%%%%%%%%%%%

\section{
%Isotopies on ${\mathcal O}$
Schur's Lemma}\label{S:Isotopies}

Let $\{B_t\}_{t\in[0,1]}$ be a family of elements in ${\frak g}$.
This family generates time-dependent vector fields on $G$,
${\mathcal O}$ and ${\mathcal F}$, which give rise to evolution
equations for several sorts of objects. In Propositions
\ref{PropstB1}, \ref{Prop2en1} and \ref{ProbfFt} we state
 properties of the solutions to these equations  that we use
 later.

By ${\varphi}_t$ we denote the isotopy on ${\mathcal O}$
determined by the time-dependent vector field $X_{B_t}$; that is,
\begin{equation}\label{isotopy'}
\frac{d\,\varphi_t}{dt}=X_{B_t}\circ\varphi_t,\;\;\;
\varphi_0={\rm id}.
\end{equation}

On the other hand the time-dependent vector field $Y_{B_t}$ on
${\mathcal F}$ defines a flow ${\mathbf H}_t$; that is, the family
of diffeomorphisms of ${\mathcal F}$  determined by
\begin{equation}\label{varphiflow'}
\frac{d\,{\bf H}_t(p)}{dt}=Y_{B_t}\big( {\bf H}_t(p) \big),\;\;\;
{\bf H}_0={\rm Id}.
\end{equation}
 %We call ${\bf H}_t$ the action integral along the isotopy
%$\varphi_t$, relative to the integral datum $\Lambda$.

Given $s\in{\mathcal B}_1$, we define a family of maps $\{s_t:
G\to H\}_{t}$  by the equations
 \begin{equation}\label{trnsportG'}
  \frac{d\,s_t}{dt}=-R_{B_t}s_t,\;\;\; s_0=s.
 \end{equation}

\begin{Prop}\label{PropstB1}
If the family $\{s_t\}$ is solution of (\ref{trnsportG'}), then
$s_t\in{\mathcal B}_1$ for all $t$.
\end{Prop}

{\it Proof.}
 \begin{equation}\label{LZst}
\frac{d}{dt}(L_Zs_t)=L_Z\big(\frac{d\,s_t}{dt}\big)=-L_ZR_{B_t}(s_t)=-
 R_{B_t}L_Z(s_t).
 \end{equation}
If $Z\in{\frak u},$ then $L_Z\,s_0=0$. The uniqueness of solutions
of the first order
 differential equation (\ref{LZst})  implies $L_Z\,s_t=0$, for all
 $t$.

On the other hand, given $l\in L$, on the space of smooth maps
$h:G\to H$ we define the operator $\alpha$ by
$$\alpha(h)=h\circ R_l-\Phi(l^{-1})h,$$
$R_l$ being the right multiplication in $G$ by $l$. If $A\in{\frak
g}$, it is straightforward to check
 %then
 % $$\alpha\big(R_A h
% \big)(g)=\frac{d}{dt}\bigg|_{t=0}\Big(h(e^{tA}gl)-\Phi(l^{-1})h(e^{tA}g)\Big)=
% \frac{d}{dt}\bigg|_{t=0}(\alpha h)(e^{tA}g)=R_A(\alpha h)(g).$$
 %That is,
 \begin{equation}\label{alpharA}
 \alpha R_A=R_A\alpha.
  \end{equation}
 If $s_t$ is
 solution of  (\ref{trnsportG'}), it follows from
 (\ref{alpharA})
 $$\frac{d}{dt}\alpha(s_t)=\frac{d\, s_t}{dt}\circ R_l-\Phi(l^{-1})\frac{d\,
 s_t}{dt}= \alpha(-R_{B_t}s_t)=-R_{B_t}(\alpha
 s_t).$$
 As $\alpha (s)=0$ since
$s\in{\mathcal B}_1$,  so we conclude $\alpha s_t=0$ for all $t$.
Thus $s_t\in {\mathcal B}_1$.

\qed

\smallskip

Given a section $\sigma\in{\mathcal B}_2$ we can consider the
family of sections $\sigma_t$ of ${\mathcal V}$ defined by the
following equations
\begin{equation}\label{fracdsigmat'}
\frac{d\,\sigma_t}{dt}={\mathcal P}_{B_t}\sigma_t,\;\;\;
\sigma_0=\sigma.
 \end{equation}
Similarly, given $f\in{\mathcal B}_3$, let $f_t$ be the set of
functions $f_t:{\mathcal F}\to H$ such that
\begin{equation}\label{definft'}
\frac{d\,{f}_t(p)}{dt}=-Y_{B_t}(p)(f_t),\;\;\;\;
  f_0=  f.
  \end{equation}
By Theorem \ref{theorempiprime} together with the preceding
Proposition one has
\begin{Prop}\label{Prop2en1}
Let $\sigma_t$ be the solution of (\ref{fracdsigmat'}) and $f_t$
the solution of (\ref{definft'}). If $\sigma$ and $f$ are
associated with $s\in{\mathcal B}_1$, then $\sigma_t$ and $f_t$
are associated with $s_t$, solution of (\ref{trnsportG'}). In
particular $\sigma_t\in{\mathcal B}_2$ and $f_t\in{\mathcal B}_3.$
\end{Prop}

%Given a smooth function $ f:{\mathcal F}\to H$, we put
Given a $f\in {\mathcal B}_3$,
  %smooth function $ f:{\mathcal F}\toH$,
  we put
\begin{equation}\label{defitionft}
\Hat f_t=f\circ {\bf H}_t^{-1}.
 \end{equation}
  We have the following proposition.
  % It is straightforward to
 %prove the following Proposition
\begin{Prop}\label{ProbfFt}
 The set of functions $\Hat f_t$ defined by (\ref{defitionft}) satisfies
 $$\frac{d\,\Hat f_t (p)}{dt}=-Y_{B_t}(p)(\Hat f_t),\;\;\;\;\Hat
f_0=f. $$
 Hence
 %$\Hat f_t=f_t$.
 $\Hat f_t$ is the solution of (\ref{definft'}).
  %Moreover, if $Z_B(f)=0$ then $Z_B(f_t)=0$.
\end{Prop}
{\it Proof.} By (\ref{varphiflow'})
 \begin{equation}\label{auxprobfFt}
 \frac{d}{du}\bigg|_{u=t}{\bf
H}_u({\bf H}_t^{-1}(p))=Y_{B_t}(p);
 \end{equation}
 that is, $Y_{B_t}(p)$ is the vector defined by the curve $\{{\bf H}_u({\bf
 H}_t^{-1}(p))
 \}_u$ at $u=t$.
 On the other hand, by (\ref{defitionft})
$$\frac{d\,\Hat f_t(p)}{dt}=W(f),$$
  where $W$ is the tangent vector to
${\mathcal F}$ at ${\bf H}_t^{-1}(p)$ defined by the curve
 $\{{\bf H}_u^{-1}(p)\}_u$. Since ${\bf H}_u({\bf H}_u^{-1}(p))=p$
 for all u, it turns out that $Y_{B_t}(p)=-({\bf H}_t)_*(W)$.
 So
 $$-Y_{B_t}(p)(\Hat f_t)=W(\Hat f_t\circ{\bf
 H}_t)=W(f)=\frac{d\,\Hat f_t(p)}{dt}.$$

\qed

\smallskip

If we integrate the family   $\{B_t\}$ we will obtain the
solutions of differential equations (\ref{isotopy'}) and
(\ref{varphiflow'}). That is, we define the curve $b_t$ in $G$ by
the conditions
 \begin{equation}\label{dotg}
 \dot b_tb_t^{-1}=B_t, \;\;\;b_0=e.
 \end{equation}
Then the isotopy $\varphi_t$ determined by (\ref{isotopy'}) is the
  multiplication by $b_t$; that is,
\begin{equation}\label{psit}
\varphi_t(g\cdot\eta)=b_t\cdot(g\cdot\eta).
\end{equation}
Analogously, the bundle diffeomorphism ${\bf H}_t$ defined in
(\ref{varphiflow'}) is the left multiplication by $b_t$ in
${\mathcal F}$,
\begin{equation}\label{bfFt}
{\bf H}_t={\mathcal L}_{b_t}.
\end{equation}

\begin{Prop}\label{Propst} The solution of (\ref{trnsportG'}) is
$s_t=s\circ L_{b_t^{-1}},$ where $b_t$ is determined by
(\ref{dotg}).
 %Moreover $s_t\in{\mathcal B}_1$
\end{Prop}

{\it Proof.}  We write $\tilde s_t:=s\circ L_{b_t^{-1}}$. Given
$g\in G$
$$\frac{d\,\tilde s_t}{dt}(g)=\frac{d}{du}\bigg|_{u=t}s(b_u^{-1}g).$$
But the curves $\{b_u^{-1}\}_u$ and $\{b_t^{-1}b_u b_t^{-1} \}_u$
in the group $G$ define opposite tangent vectors  at $u=t$.
 So
 $$\frac{d\,\tilde s_t}{dt}(g)=-\frac{d}{du}\bigg|_{u=t}s(b_t^{-1}b_u
 b_t^{-1}g).$$
 On the other hand,
 $$R_{A_t}(\tilde s_t)(g)=\frac{d}{du}\bigg|_{u=t}\tilde s_t(b_ub_t^{-1}g)=
\frac{d}{du}\bigg|_{u=t}s(b_t^{-1}b_u
 b_t^{-1}g).$$
  That is, $s\circ L_{b_t^{-1}}$ satisfies
 (\ref{trnsportG'}).
  \qed

\medskip

{\bf Proof of Theorem \ref{kappatransport}.}
 Now $g_1\in Z(G)$ and
the isotopy $\psi_t$ defined in (\ref{isotopy}) is closed; that
is, $\psi_1={\rm id}$.
 %Furthermore, as $\pi$  is an unitary irreducible
%representation of $G$ on Hilbert space completion of  ${\mathcal
%B}_1$.
 By Proposition \ref{Propst}, (\ref{equivariantbis}) and
(\ref{kappadefinition}) one has
$$s_1=s\circ L_{{g_1}^{-1}}=\pi_1(g_1)(s)=\kappa s,$$
for any $s\in{\mathcal B}_1$, which proves item (a).

By (\ref{sigma-s}) and Proposition \ref{Prop2en1},
  the result stated in (a)   expressed in terms
of the solutions to
 (\ref{fracdsigmat}) gives (b).

Moreover
$$\kappa s(g)=s_1(g)=s(g_1^{-1}g)=s(gg_1^{-1})=\Phi(g_1)s(g);$$
 that is,
\begin{equation}\label{Phig1}
 \Phi(g_1)=
 \kappa\,{\rm Id}.
  \end{equation}

% Denoting by ${\bf F}_t$ the flow in ${\mathcal F}$ defined by
 %$Y_{A_t}$,
 It follows from (\ref{bfFt}) and (\ref{Phig1}) that
    \begin{equation}\label{boldF1}{\bf F}_1[g,\,\alpha]=[g_1g,\alpha]
    =[gg_1,\alpha]=[g,\,\Phi(g_1)\alpha]=\kappa[g,\alpha],
 \end{equation}
and (d) is proved.

 From Proposition \ref{ProbfFt} together with
(\ref{bfFt}) and (\ref{boldF1}), it follows
$$f_1[g,\,\alpha]=(f\circ{\bf F}_1^{-1})[g,\,\alpha]=f[g,\,\kappa^{-1}\alpha]=\kappa
f[g,\,\alpha],$$
 which  proves (c).
  \qed

\medskip

  Let us assume that there is a fixed point
$x_0\in{\mathcal O}$ for the isotopy $\{\psi_t\}$ defined in
(\ref{isotopy}); that is, $\psi_t(x_0)=x_0$, for all $t$. So
$X_{A_t}(x_0)=0$ and (\ref{fracdsigmat}) evaluated at $x_0$
reduces to
$$\frac{d\,\sigma_t(x_0)}{dt}=h_{A_t}(x_0)\sigma_t(x_0),\;\;\;\,
\sigma_0(x_0)=\sigma(x_0).$$ This is a differential linear
equation for $v_t:=\sigma_t(x_0)\in H.$ Let $M(t)\in\frak {gl}(H)$
be the ``fundamental matrix" of this linear equation,
 in other words
  $$\frac{d\, M(t)}{dt}=h_{A_t}(x_0)M(t),\;\;\; M(0)={\rm Id}.$$
 By Theorem \ref{kappatransport} it follows
\begin{equation}\label{M(1)=}
 M(1)=\kappa\,{\rm Id}.
 \end{equation}

\smallskip

 {\bf Corollary to Theorem \ref{kappatransport}.}
{\it If $A_t=A$ for all $t$, and $x_0$ is a fixed point of the
isotopy $\{\psi_t\}_t$, then
 $$\kappa\,{\rm Id}= {\rm exp}\big( h_A(x_0) \big).$$}

%%%%%%%%%%%%%%%%%%%%%%%%%%%%%%%%%%%%%%%%%%%%%%%%%%%%%%%%%%%%%%%%%%%%%%%%%%%%%%%%%%%%%%%%%%%%%%%%%%%%%%%%%%%%%%%%%%%%%%%%%%
%%%%%%%%%%%%%%%%%%%%%%%%%%%%%%%%%%%%%%%%%%%%%%%%%%%%%%%%%%%%%%%%%%%%%%%%%%%%%%%%%%%%%%%%%%%%%%%%%%%%%%%%%%%%%%%%%%%%%%%%%%%%%%

\smallskip

\subsection{ Case when  ${\;\rm dim}\,H=1$.}\label{SubsectionCase}

Now the bracket and the adjoint action  in $\frak {gl}(H)$ are
trivial. It follows from (\ref{Curvature})
$${\bf K}_{[g,\alpha]}(Y_B,\,Y_C)=-\phi(g^{-1}[B,\,C])=-h_{[B,\,C]}(g),$$
and ${\bf K}$ projects a closed $2$-form $K$ on ${\mathcal O}.$ We
denote by $\omega:=-K$; that is,
 \begin{equation}\label{omega}
 \omega(X_A,X_B)(g\cdot\eta)=h_{[A,\,B]}(g).
  \end{equation}
In this case   (\ref{DHathB}) reduces to
 \begin{equation}\label{dhb=}
 dh_B=\iota_{X_B}\omega.
 \end{equation}

  In this context $\psi_t$ defined in (\ref{isotopy})
 is an isotopy which  determines the time-dependent Hamiltonian
 $h_{A_t}$ through the form $\omega$ \cite{Mc-S}.

  \medskip

{\bf Proof of Theorem \ref{integralaccion}.}
  Let $\mu$ be a local
frame for the line bundle ${\mathcal V}$. The solution $\sigma_t$
to (\ref{fracdsigmat}) can be written $\sigma_t=m_t \mu$, where
$m_t$ is a complex function defined on an open set of ${\mathcal
O}$. Then (\ref{fracdsigmat}) gives rise to
\begin{equation}\label{dmt}
\frac{d\,m_t}{dt}=-\gamma(X_{A_t})m_t-X_{A_t}(m_t)+h_{A_t}m_t,
\end{equation}
where $\gamma$ is the connection form of ${\mathcal V}$ in the
frame $\mu$.
 Given  an arbitrary point  $q$ of ${\mathcal O}$, then
 $\{q_t:=\psi_t(q)\}_t$ is a  closed curve on ${\mathcal O}$. If
 $q$ belongs to the domain of $\mu$, we define $m'_t:=m_t(q_t)$.
 If we evaluate  (\ref{dmt}) at the point $q_t$ we obtain
 $$\frac{d\,m'_t}{dt}=\big(-\gamma_{q_t}(X_{A_t})+h_{A_t}(q_t)
 \big)m'_t.$$
So
$$m'_t=m'_0\,{\rm exp}\Big(\int_0^t\big(-\gamma_{q(u)}(X_{A_u})+h_{A_u}(q_u)  \big)du
\Big).$$

On the other hand, we can consider on ${\mathcal O}$ the Kirillov
symplectic structure \cite{Kir0}, then $\psi_t$ is a Hamiltonian
isotopy with respect to this structure,
%If $\psi_1={\rm id}$, then
%$\{\psi_t\}_t $ is a loop in corresponding Hamiltonian group,
 and
consequently the evaluation closed curve $\{ q_t\}$ is
nullhomologous (Lemma 10.31 in \cite{Mc-S}), that is, it is the
boundary of a $2$-chain.
%Therefore, if $g_1\in Z(G)$, by
By Stokes' theorem
$$m'_1=m'_0\,{\rm
exp}\Big(\int_S\omega+\int_0^1 h_{A_t}(q_t)dt \Big),$$
 where $S$ a $2$-chain whose boundary is the curve is $\{q_t\}_t$.
 By Theorem \ref{kappatransport}
 \begin{equation}\label{kappaintegralaction}
\kappa= {\rm exp}\Big(\int_S\omega+\int_0^1 h_{A_t}(q_t)dt \Big).
\end{equation}
  \qed

\smallskip

  {\it Remarks.}
   %From (\ref{kappaintegralaction}) one deduces that the righthand side of this equality
 %is independent of the point $q$.
  %This independence is   known \cite{V},  and
 %the value of this exponential is called the action integral around the closed isotopy $\psi_t$
 %(see \cite{Mc-S}, \cite{We}). In view this fact and item (d) of
 %Theorem \ref{kappatransport}, the flow
 %${\bf H}_t$ (defined in (\ref{varphiflow'})) may be considered as the
 %``action integral" along the isotopy $\varphi_t$,
 %(defined in (\ref{isotopy'}))
 %relative to the integral datum $\Lambda$.
If $q$ is a fixed point for $\psi_t$ and $A_t=A$ for all $t$, from
Theorem \ref{integralaccion} it follows  $\kappa={\rm
exp}(h_A(q))$; this agrees with
 Corollary to Theorem \ref{kappatransport}.
%with Corollary \ref{Corkappa}.
 %(\ref{M(1)At=A}).

 By (\ref{kappadefinition}) the exponential in the statement of
Theorem \ref{integralaccion} depends only on the final point of
the curve $g_t$ and it is independent of the family $A_t$ defined
by $g_t$.

%%%%%%%%%%%%%%%%%%%%%%%%%%%%%%%%%%%%%%%%%%%%%%%%%%%%%%%%%%%%%%%%%%%%%%%%%%%%%%%%%%%%%%%%%%%%%%%%%%
%%%%%%%%%%%%%%%%%%%%%%%%%%%%%%%%%%%%%%%%%%%%%%%%%%%%%%%%%%%%%%%%%%%%%%%%%%%%%%%%%%%%%%%%%%%%%%%%%%%

\subsection{The   character}\label{Subsection:the character}

A slight modification of the preceding developments allows us to
prove the formula for the character   given in Theorem
\ref{theoremcharacter}.

\smallskip

{\bf Proof of Theorem \ref{theoremcharacter}.}   Let $c_t$
 %Let $\tilde g$ be an element whose conjugacy class meets $L$; that
 %is, such that $a^{-1}\tilde g a\in L$ with  $a\in G$. By $c_t$ we
denote a path in $L$ with $c_1=a^{-1}\tilde g a$. Then
$g_t:=ac_ta^{-1}$ is a path in $G$ with $g_1=\tilde g$. The point
$x_0=a\cdot\eta\in{\mathcal O}$ is a fixed point for the isotopy
on ${\mathcal O}$
 defined by multiplication by $g_t$. We put $A_t\in{\frak g}$ and
  $C_t\in{\frak l}$ for the velocity paths associated with $g_t$
  and $c_t$, respectively. So $A_t={\rm Ad}_aC_t$.

 We denote by $V_2$ the subspace of ${\mathcal
  B}_2$ on which $\tau$ is defined.
  The action of $\tau(\tilde g)$ on a section $\sigma\in{V}_2$ is the section
   $\sigma(1)$ determined by the equations
  $$\frac{d\,\sigma(t)}{dt}={\mathcal P}_{A_t}(\sigma(t)),\;\;\;
  \sigma(0)=\sigma.$$
  Evaluating these equations at the point $x_0$, and taking into
  account that $X_{A_t}(x_0)=0$, one obtains
 \begin{equation}\label{dsigma(t)}
  \frac{d\,\sigma(t)(x_0)}{dt}=
h_{A_t}(x_0)\big(\sigma(t)(x_0)\big),\;\;\;
  \sigma(0)(x_0)=\sigma(x_0).
   \end{equation}
  As ${\dim}\, H=1$, it follows from (\ref{dsigma(t)}) that
  \begin{equation}\label{sigma(1)}
  \sigma(1)(x_0)={\rm exp}\Big(\int_0^1h_{A_t}(x_0)
  dt\Big)\,\sigma(x_0),
  \end{equation}
for any $\sigma\in{V}_2$.

If ${\rm dim}\,\tau=m$, let $\sigma_1,\dots,\sigma_m$ be a basis
of ${V}_2$, then
$$\sigma_j(1)=\tau(\tilde g)\sigma_j=\sum_iM_{ij}\sigma_i,$$
with $M_{ij}\in{\mathbb C}.$ By (\ref{sigma(1)})
$$M_{ji}=\delta_{ji}\,{\rm exp}\Big(\int_0^1h_{A_t}(x_0)
  dt\Big),$$
and the proof is complete.
 \qed

\smallskip

{\it Remark.} As $h_{A_t}(x_0)=\phi(a^{-1}\cdot A_t)=\phi(C_t)$
the formula for the character can be written
$$\chi_{\tau}(\tilde g)=m\,{\rm
exp}\Big(\phi\big(\int_0^1 C_t dt\big)\Big).$$

 If $\tilde c:=a^{-1}\tilde g a$ equals $e^C$, with $C\in{\frak l}$, then we
 can take $C_t=C$ for all $t$ and
 \begin{equation}\label{chireduced}
  \chi_{\tau}(\tilde g)=m\,\Phi(\tilde c).
 \end{equation}

  (\ref{chireduced}) can  also be deduced by considering
 $\tau$ as a representation on a subspace $V_1$ of ${\mathcal B}_1$. Given $s\in{V}_1$,
 $$(\tau(\tilde g)s)(a)=s(a\tilde c^{-1})=\Phi(\tilde c)s(a).$$
  From this formula it follows (\ref{chireduced}).

%%%%%%%%%%%%%%%%%%%%%%%%%%%%%%%%%%%%%%%%%%%%%%%%%%%%%%%%%%%%%%%%%%%%%%%%%%%%%%%%%%%%%%%%%%%%%%%%%%%%%%%%%%%%%%%%%
%%%%%%%%%%%%%%%%%%%%%%%%%%%%%%%%%%%%%%%%%%%%%%%%%%%%%%%%%%%%%%%%%%%%%%%%%%%%%%%%%%%%%%%%%%%%%%%%%%%%%%%%%%%%%%%%%

\subsection{The infinitesimal character}\label{Subsection:infinitesimal}

% The center of
%the universal enveloping algebra $U({\frak g_{\mathbb C}})$ of
%${\frak g_{\mathbb C}}$
 % will be denoted by ${\mathcal
 % Z}({\frak g_{\mathbb C}})$.
% Let
%${\frak h}$ be any Cartan subalgebra of ${\frak g}$, By
%Harish-Chandra isomorphism ${\mathcal Z}({\frak g})$ can be
%identified with the elements of the symmetric algebra $S({\frak
%h})$ invariant under the Weyl group \cite{Kn-Vo}.

 Now we consider the ``representation" of the
associative algebra $U({\frak g}_{\mathbb C})$ induced by $\pi'$
on the space $({\mathcal B}_2)_K$ of $K$-finite vectors in
${\mathcal B}_1$. Since $\pi$ is unitary and irreducible, each
element of ${\mathcal Z}({\frak g}_{\mathbb C})$ acts as a scalar
operator
 %in the representation $\pi^{U({\frak g})}$
 (see \cite{kn0}, Corollary 8.14).
 By
$\chi:{\mathcal Z}({\frak g}_{\mathbb C})\to{\mathbb C}$ we denote
the corresponding infinitesimal character.

Let ${\frak h}$ be a Cartan subalgebra of ${\frak g}_{\mathbb C}$
contained in ${\frak l}_{\mathbb C}$. We denote by $\Delta^+$ a
set of positive roots for the pair $({\frak g}_{\mathbb
C},\,{\frak h})$. For $\alpha\in\Delta^+$, $E_{\alpha}$ will be a
basis for the corresponding root space. According to Lemma 8.17 in
\cite{kn0}, if $Z\in{\mathcal Z}({\frak g}_{\mathbb C})$ then
$Z\in U({\frak h})\oplus {\mathcal P},$ where
 $${\mathcal P}=\sum_{\alpha\in \Delta^+}U({\frak g}_{\mathbb
 C})E_{\alpha}.$$
 The projection of $Z$ into $U({\frak h})$ will be denoted $\Hat Z$.

Let $V\subset({\mathcal B}_2)_K$ be an irreducible representation
of the maximal compact subgroup $K$ which occurs in $\pi_{|K}$.
Now we consider $s_0\in V$ a highest weight vector of the
representation $V$, so $\pi'({E_{\alpha}})s_0=0$ for all
$\alpha\in\Delta^+$. Thus, if $Z\in{\mathcal Z}({\frak g}_{\mathbb
C})\cap{\mathcal P}$ then the action of $Z$ on $s_0$ vanishes. As
$\pi'(Z) s_0=\chi(Z) s_0$, it follows $\chi(Z)=0$. We have the
following proposition
\begin{Prop}\label{PropchiP}
If $Z\in{\mathcal Z}({\frak g}_{\mathbb C})\cap{\mathcal P}$, then
$\chi(Z)=0.$
\end{Prop}
   To prove Theorem
\ref{7rule} we need the following Lemma
 \begin{Lem}\label{LemRALB}
  With $\phi$ denoting the extension of the map (\ref{definphi})  to
${\frak g}_{\mathbb C}$ and $1$ the identity  element of $G$, we
have

  (i) If $Y\in{\frak h}$ and $s\in{\mathcal B}_1$, then
 $\pi'(Y)\,s=\phi(Y)\,s.$

(ii) If $Y,W\in{\frak h}$, then
 \begin{equation}\notag
     %(R_AR_B\,s)(1)=\phi(A)\phi(B)s(1).
     (\pi'(Y)\pi'(W)\, s)(1)=\phi(Y)\phi(W)\,s(1)
 \end{equation}

 \end{Lem}

{\it Proof.} Since ${\frak h}\subset {\frak l}_{\mathbb C}$ we can
assume that $Y,W$ are elements of ${\frak l}$. The item (i)
follows from  (\ref{piprime}) together with the fact that $s$ is
$\Phi$-equivariant.

\smallskip

 If $Y,W\in {\frak l}$, then
 $s(e^{uW}e^{tY})=\Phi(e^{-tY})\Phi(e^{-uW})s(1)$. Hence
 $$(R_YR_W\,s)(1)=\frac{d}{dt}\bigg|_{t=0}\frac{d}{du}\bigg|_{u=0}\Phi(e^{-tY})\Phi(e^{-uW})\,s(1)=
 \phi(Y)\phi(W)\,s(1),$$
and (ii) follows.
 \qed

 %Note that $\phi$ is a Lie algebra homomorphism which satisfies
%(\ref{invphi}), so $\phi(A)\phi(B)=\phi(B)\phi(A)$, for
 %$A,B\in{\frak l}_{\mathbb C}$.

\smallskip

Given $\{Y_1,\dots, Y_r\}$ a basis of ${\frak h}$, and
$Z\in{\mathcal Z}({\frak g}_{\mathbb C})$, then there exist a
polynomial $q(Y_1,\dots, Y_r)$ as in (\ref{Polynomial}) such that
$\Hat Z=q(Y_1,\dots Y_r)$.

\begin{Prop}\label{Propchi}
With the above notations
$$q(\phi(Y_1),\dots,\phi(Y_r))=\chi(Z)\,{\rm Id}.$$
\end{Prop}

 {\it Proof.} By
  %Lemma \ref{LemRALB} and
 Proposition \ref{PropchiP}, the operator $\pi'(Z)$ associated
to $Z$ is
 $q(-R_{Y_1},\dots,-R_{Y_r})$. By Lemma \ref{LemRALB}, if
   $s\in{\mathcal B}_1$ then
 $$\big(q(-R_{Y_1},\dots,-R_{Y_r})\,s\big)(1)=q\big(\phi(Y_1),\dots,\phi(Y_r)
 \big)\,s(1).$$
 As $q(-R_{Y_1},\dots,-R_{Y_r})=\chi(Z)\,{\rm Id}$, we obtain the
 proposition.
 \qed

 \begin{Thm}\label{Theoremcharacter0}
Given $Z\in{\mathcal Z}({\frak g}_{\mathbb C})$, if $\Hat
Z=q(Y_1,\dots Y_r)$ and $h_k:=h_{Y_k}$, the Hamiltonian map
associated with $Y_k$,  then the function
$$q(h_1,\dots, h_r):{\mathcal O}_{\eta}\to \frak{gl}(H)$$
  takes at the point
$\eta$ the value $\chi(Z)\,{\rm Id}$.
\end{Thm}
{\it Proof.} It follows from Proposition \ref{Propchi} and
(\ref{definh}).

\qed

\smallskip

 {\bf Proof of Theorem \ref{7rule}.} If $\tilde q$ is any
  polynomial in the variables $Y_1,\dots,Y_r$, from
(\ref{lifhA}) one deduces that
 \begin{equation}\label{homogeneous}
 \tilde q({\bf
h}_1,\dots,{\bf h}_r)[g,\,\alpha]=\alpha^{-1}\big(\tilde
q(h_1,\dots,h_r)(g\cdot\eta)\big)\alpha.
\end{equation}

 %If $q$
 %defines the element $\Hat Z$,  by Theorem \ref{Theoremcharacter0}
On the other hand, $q(h_1,\dots,h_r)(\eta)$ is a multiple of
identity, by Theorem \ref{Theoremcharacter0}. From this fact
together with (\ref{homogeneous}) we   deduce
 $$ q({\bf h}_1,\dots,{\bf
h}_r)[1,\,\alpha]=q(h_1,\dots,h_r)(\eta),$$
 for any
$[1,\,\alpha ]\in{\mathcal F}.$ The theorem follows from Theorem
\ref{Theoremcharacter0}. \qed

\smallskip

Let $x_0=g\cdot\eta$ be a point of ${\mathcal O}$ and let ${\frak
h}'$ be a Cartan subalgebra of ${\frak g}_{\mathbb C}$ such that,
$g^{-1}\cdot {\frak h}'\subset {\frak l}_{\mathbb C}$. A
generalization of Theorem \ref{7rule} is the following proposition

\begin{Prop}\label{Propultima}
If $Y'_1,\dots,Y'_r$ is a basis of ${\frak h}'$ and the polynomial
$q'(Y'_1,\dots,Y'_r)$ is the projection of $Z\in{\mathcal
Z}({\frak g}_{\mathbb C})$ into $U({\frak h}')$, then the function
$q'({\bf h}_{Y'_1},\dots,{\bf h}_{Y'_r})$ is constant on the fiber
of ${\mathcal F}$ over $x_0$, and its value on this fiber is
$\chi(Z){\rm Id}$.
 \end{Prop}

%%%%%%%%%%%%%%%%%%%%%%%%%%%%%%%%%%%%%%%%%%%%%%%%%%%%%%%%%%%%%%%%%%%%%%%%%%%%%%%%%%%%%%%%%%%%%%%%%%%%%%%%%%%%%%%%%%%%%%%%%%%
%%%%%%%%%%%%%%%%%%%%%%%%%%%%%%%%%%%%%%%%%%%%%%%%%%%%%%%%%%%%%%%%%%%%%%%%%%%%%%%%%%%%%%%%%%%%%%%%%%%%%%%%%%%%%%%%%%%%%%%%

\medskip

\section{Physical interpretations}\label{Sect.Physical}

As it is well known Geometric Quantization \cite{WO} is a
mathematical procedure for understanding the relation between a
classical physical system and its ``quantization". From the
mathematical point of view the classical phase space is a
symplectic manifold $(M,\alpha)$, and the set of rays of a Hilbert
space ${\mathcal H}$ is the mathematical model for the space of
states of the quantum system. The manifold $(M,\alpha)$ is said to
be quantizable if the cohomology class of $\alpha/(2\pi)$ is
integral. In this case there exists a Hermitian line bundle
${\mathcal L}$ on $M$  equipped with a connection whose curvature
is $-i\alpha$.
 ${\mathcal L}$ is called a ``prequantum bundle".
%The next step in the
For the construction of the Hilbert space ${\mathcal H}$ from the
quantizable manifold $M$ one  fixes a polarization ${\frak P}$ on
$M$, then ${\mathcal H}$ is a subset of
  the  space of sections of ${\mathcal L}$
polarized with respect to ${\frak P}$ (see \cite{WO} and \cite{Sn}
for the details omitted in this schematic summary).

The coadjoint orbit ${\mathcal O}$ of $\eta\in{\mathfrak g}^*$
supports a canonical symplectic structure $\Hat\omega$, the
Kirillov form \cite{Kir0}. Denoting by $L$   the stabilizer of
$\eta$, the orbit  ${\mathcal O}$  admits a $G$-invariant
prequantization iff the linear operator $i\eta:{\mathfrak l}={\rm
Lie}\,(L)\to i{\mathbb R}$ is integral, in the sense that there
exists a character $\Lambda:L\to U(1)$ whose derivative is $i\eta$
\cite{Kos}.
 In this case the
corresponding prequantum bundle is ${\mathcal
L}=G\times_{\Lambda}{\mathbb C}$. Since   the group $G$ acts by
translation on the orbit, it is reasonable to impose the
quantization to have a $G$-invariant Hilbert space structure.
  In general it is not possible to integrate the absolute value of sections of ${\mathcal L}$
in a translation-invariant way, since ${\mathcal O}=G/L$ does not
admit a measure invariant under the action of $G$
 (see p. 537 \cite{Kn}). To define such an   integration it is necessary  to
consider a prequantum bundle  different from ${\mathcal L}$;
specifically, one takes the bundle $\mathcal V$ determined
 by the character
$\Phi=\Lambda\cdot\Delta$, where $\Delta^2$ is the modular
function on $G/L$. If $\sigma_1,\sigma_2$ are compactly supported
sections of $\mathcal V$, and $s_1,s_2:G\to{\mathbb C}$ are the
corresponding $\Phi$-equivariant functions, then
$m(g):=s_1(g)\overline{s_2(g)}$ satisfies
\begin{equation}\label{DeltaEquivariant}
m(gl)=\Delta^2(l^{-1})m(g),\;\;\;\text{for\; all}\;l\in L.
\end{equation}
(That is, $m$ defines a section    of the bundle  of densities on
$G/L$.) The functions on $G$ which satisfies
(\ref{DeltaEquivariant}) have a translation invariant integral
``over G/L" (see p. 65 \cite{Vo0}, p. 41 \cite{K-T}). Then
$\langle \sigma_1,\,\sigma_2 \rangle:=\int_{G/L}m$ defines a
$G$-invariant product of compactly supported sections of $\mathcal
V$.

 On the other hand, if $\eta$ is a
hyperbolic element then the orbit ${\mathcal O}$ possesses the
polarization determined by the subalgebra ${\frak u}$ defined in
(\ref{fraku}). So our space ${\mathcal B}_2$, defined in
(\ref{B2}), is  a $G$-invariant quantization of the orbit. By
Proposition \ref{PropBi} the spaces ${\mathcal B}_j$, $j=1,2,3$,
can be considered as equivalent $G$-invariant quantizations  of
${\mathcal O}$.

\smallskip

The Kirillov  symplectic form $\Hat\omega$ is defined by
$$\Hat\omega_{g\cdot\eta}(X_A,\,X_B)=\eta(g^{-1}\cdot[A,B]),$$
and the Hamiltonian function associated to $A$ is $\Hat
h_A(g\cdot\eta)=\eta(g^{-1}\cdot A)$. From (\ref{definh}),
(\ref{omega}) and (\ref{definphi}) one obtains
$$\omega=i\Hat\omega+\tilde\omega,\;\;\;\; h_A=i\Hat h_A+\tilde h_A, $$
where
$$\tilde\omega_{g\cdot\eta}(X_A,\,X_B)=\delta(g^{-1}\cdot[A,B]),\;\;\;\;
\tilde h_A(g\cdot\eta)=\delta(g^{-1}\cdot A).$$
 $\,\tilde\omega$ is not a
 symplectic form (because it is degenerate), but the analogous  relations to (\ref{dhb=})   with tildes and
 with hats are also valid.

 %Moreover, the analogous  relations to (\ref{dhb=})   with tildes and
 %with hats are also valid. (However $\tilde\omega$ is not a
 %symplectic form, because it is degenerate).

 \smallskip

 The cotangent bundle $M=T^*P$ to a manifold carries a canonical
 $1$-form $\beta_0$ \cite{Mc-S}, and $\alpha_0:=-d\beta_0$ defines a
 symplectic structure on $M$. If $q:t\in[0,\,1]\to M$ is a curve
 and $h_t:M\to {\mathbb R}$ is a time dependent Hamiltonian on
 $M$, the action integral along the curve $q(t)$ is defined by the
 following formula \cite{A-M}, \cite{Mc-S}
 $$\int_0^1\big(-\beta_0(\Dot q(t))+h_t(q(t))  \big)dt.$$

For a general symplectic manifold $(M,\,\alpha)$ the time
dependent Hamiltonian $ h_t$, with $t\in [0,1]$, determines a time
dependent Hamiltonian vector field $X_t$, which in turn defines an
isotopy of symplectomorphisms  $\xi_t$. If $\xi_1={\rm id}$ (that
is, $\{\xi_t\}_{t\in[0,1]}$ is a loop in ${\rm Ham}(M)$, the
Hamiltonian group of $M$ \cite{Mc-S}), then evaluation curve
$\{\xi_t(p)\}_t$ is nullhomotopic, for all point $p\in M$
\cite{L-M-P}. Hence the action integral around this curve can be
written
\begin{equation}\label{ActionInte}
 \Hat{\mathcal A}(\xi):=\int_S\alpha+\int_0^1h_t(\xi_t(p))dt,
\end{equation} $S$ being a $2$-chain whose boundary is the curve
$\{\xi_t(p)\}_t$. It is known that the value of (\ref{ActionInte})
is independent of the point $p$ \cite{V}.

In the case when the manifold is a coadjoint orbit ${\mathcal O}$
and the loop $\{\psi_t\}$ in  ${\rm Ham}({\mathcal O})$ is defined
as in (\ref{isotopy}),
% that is, by action of the group $G$,
 one can also consider the ``action integral" $\tilde{\mathcal A}(\psi)$
%defined as in (\ref{ActionInte}) but
 %using
 defined by means of the $2$-form
$\tilde\omega$ and the ``Hamiltonian"
 $\tilde h_{A_t}$
 %(we use quotation marks  because  $\tilde\omega$ is
 %$not a symplectic form, it is degenerate)
$$\tilde{\mathcal A}(\psi):=\int_S\tilde\omega+\int_0^1\tilde h_{A_t}(\psi_t(p))dt.$$
Thus the result stated in Theorem \ref{integralaccion} can be
written as $\kappa={\rm exp}\big(i\Hat{\mathcal
A}(\psi)\big)\times
 {\rm exp}\big( \tilde{\mathcal A}(\psi)\big).$
 %where $\tilde{\mathcal A}$ is defined as in (\ref{ActionInte})
 %but using the $2$-form $\tilde\omega$ and the Hamiltonian
 %$\tilde h_{A_t}$.
  %Obviously the invariant $\kappa$ is independent of
 %the point $p$, so $\tilde{\mathcal A}(\psi)$  does not depend on
 %the point $p$.
 Since the representation $\pi$ is unitary,
$\kappa\in U(1)$. So
 $$\kappa={\rm exp}\big(i\Hat{\mathcal A}(\psi)\big).$$
That is, the invariant of the representation $\pi$ associated to
$g_1\in Z(G)$,  by Schur's lemma, equals the exponential of $i$
times the  action integral around the loop in ${\rm Ham
}({\mathcal O})$ generated by any path in $G$ with endpoint at
$g_1$.

 In view of item (d)  of
 Theorem \ref{kappatransport} and  the above facts, we may   consider
 the flow ${\bf H}_t$ (defined in (\ref{varphiflow'}))
 as a generalized
  action integral along the isotopy $\varphi_t$,
 %(defined in (\ref{isotopy'}))
 relative to the integral datum $\Lambda$.

\smallskip

The Berry phase is a general phenomenon which may appear when a
quantum system undergoes a cyclic evolution.
 We summarize the
geometric definition of Berry phase given in \cite{We1}, where
general references can be find. Let
$\{(N_t,\epsilon_t)\}_{t\in[0,1]}$ be a family of weighted
Lagrangian submanifolds of $(M,\alpha)$ ($\epsilon_t$ being a
smooth density on $N_t$), obtained from $N=N_0$ by a set $\xi_t$
of symplectomorphisms of $M$ determined by a time dependent
Hamiltonian $h_t$. Let us assume that $(M,\alpha)$  is quantizable
and ${\mathcal L}$ is a prequantum bundle.  ${\mathcal
L}^{\times}= {\mathcal L}\setminus\{\rm zero\; section\}$ is the
corresponding principal bundle. We denote by $F_t$ the flow on
${\mathcal L}^{\times}$ generated by the vector field
$X_t^{\sharp}-W_{h_t}$, where $X_t^{\sharp}$ is the horizontal
lift of the respective Hamiltonian  vector field $X_t$, and
$W_{h_t}$ is $h_t$ times the fundamental vertical vector field on
${\mathcal L}^{\times}$. If $\sigma$ is a section of ${\mathcal
L}^{\times}_{|N}$, then $F_1(\sigma(N))$ differs from $\sigma(N)$
by  a phase $\theta$. If the Hamiltonian $h_t$ are normalized so
that $\int_{N_t}h_t\epsilon_t=0$, then $\theta$ is the Berry phase
of the loop $\{(N_t,\epsilon_t)\}_{t\in[0,1]}$ (p.142 \cite{We1}).

  By (\ref{Xsharp}) the statement (d) in Theorem
\ref{kappatransport}, when ${\rm dim}\,H=1$, can be interpreted by
saying that the invariant $\kappa$ is the Berry phase of any loop
$\{(N_t,\epsilon_t)\}_{t\in[0,1]}$ generated by the Hamiltonian
functions $h_{A_t}$, where $A_t$ is the velocity curve of any path
in $G$ with endpoint at $g_1$, and   $h_{A_t}$ is given by
(\ref{definh}).

\medskip

{\it Example.} Let $G$ be the restricted Lorentz group
$SO^+(1,\,3)$. A basis for the Lie algebra $\frak {so}^+(1,\,3)$
is $X_1,\dots,X_6$, where $X_1,X_2,X_3$ are generators of the
boosts along the axes, and $X_4,X_5,X_6$ are the generators of the
rotations around those axes. The matrix of Killing metric in the
basis $X_i$ is
$$\big(g_{ij}\big)=\big(\text{Tr}\,(\, \text{ad}\,X_i\circ \text{ad}\,X_j
)\big)
  =\,\text{diag}\,(1,1,1,-1,-1,-1),$$
   and the Casimir
element $C$ of $U({\frak g}_{\mathbb C})$ is $C=\sum
g_{ij}X^iY^j$, where $X^i=\sum _kg^{ia}X_k$ \cite{kn0}. That is,
 \begin{equation}\label{Casimir}
 C=\frac{1}{4}\Big( \sum_{1}^3X_i^2 - \sum_{4}^6X_i^2 \Big).
 \end{equation}

Let $Y=(Y_{ab})$ be a matrix in $\frak {so}^+(1,\,3)$, with
$a,b=0,1,2, 3$, and let $\eta$ denote the element in $\mathfrak
{so}^+(1,\,3)^*$ defined by $\eta(Y)=kY_{01}$, with $k\in{\mathbb
R}\setminus\{0\}$ (equivalently $\eta(X_j)=k\delta_{1j}$). The
matrix associated with $\eta$ according to (\ref{X0}) is precisely
$kX_1$, which has real eigenvalues; i.e. the coadjoint orbit
${\mathcal O}$ of ${\eta}$ is  hyperbolic. Furthermore
 ${\mathcal O}$ is $G/L$, with
$L=SO^{+}(1,\,1)\times SO(2)$. If $(A,B)\in L$, then $A$ will have
the form $A=\text{exp}\,aX_1$ and $\eta$ can be extended to a
character $\Lambda$ on $L$ by putting $\Lambda(A,B)=e^{i ka}$.

One has the following relations
$$[X_1,X_2]=X_6,\;\, [X_1,X_3]=X_5,\;\, [X_1,X_4]=0,\,
\; [X_1,X_5]=X_3,\;\, [X_1,X_6]=X_2.$$ So a basis for  the
subalgebra ${\frak u}$ defined in (\ref{fraku}) is
$X_2+X_6,\,X_3+X_5,$ and for the operators $\delta$ and $\phi$
introduced in Section 2 we have
 \begin{equation}\label{deltaX_1}
  \delta(X_1)=\text{Tr}\,(\text{ad}(X_1)_{|{\frak u}})=2, \;\;\;\;
\phi(X_1)=ik+2.
 \end{equation}
 Analogously
 \begin{equation}\label{DeltaX4}
  \delta(X_4)=0,\;\;\;\; \phi(X_4)=0.
   \end{equation}

  ${\frak h}={\frak l}_{\mathbb C}$ is a Cartan
subalgebra of ${\frak g}_{\mathbb C}$, and it is generated by
$X_1,X_4$. By (\ref{Casimir}), the projection $\Hat C$ of $C$ on
$U({\frak h})$ is $\Hat C=\frac{1}{4}(X_1^2-X_4^2).$  By Theorem
\ref{7rule}, it follows from (\ref{deltaX_1}), (\ref{DeltaX4}) and
(\ref{definh}) that the value $\chi(C)$ of the infinitesimal
character of the representation $\pi'$ associated with the orbit
${\mathcal O}$ is $(1/4)(ik+2)^2$.

For $A\in{\frak g}$, the operator  ${\mathcal P}_A=-D_{X_A}+h_A$
acting on polarized sections of the prequantum   bundle
 is the ``quantization" of the vector field on
${\mathcal O}$ determined by $A$ \cite{Sn}.
 %Since the Casimir element
 %$C$ belongs to ${\mathcal Z}({\frak g}_{\mathbb C})$,
 From the above result it turns out that the operator
$$\frac{1}{4}\Big(\sum_{i=1}^3  ({\mathcal P_{X_i}})^2 -  \sum_{i=4}^6  ({\mathcal P_{X_i}})^2  \Big),$$
associated with the Casimir element $C\in{\mathcal Z}({\frak
g}_{\mathbb C})$ is simply the multiplication by the constant
$(1/4)(ik+2)^2$.

%%%%%%%%%%%%%%%%%%%%%%%%%%%%%%%%%%%%%%%%%%%%%%
%%%%%%%%%%%%%%%%%%%%%%%%%%%%%%%%%%%%%%%%%%%%%%%%

\end{document}